\documentclass[twocolumn]{autart}

\usepackage{mathtools}
\usepackage{cuted}
\usepackage{tikz}
\usepackage{float}
\usetikzlibrary{arrows.meta}
\tikzset{>={Latex[width=1.75mm,length=1.75mm]}}
\usepackage{chngcntr}
\usepackage{apptools}
\AtAppendix{\counterwithin{lma}{section}}
\usepackage{cite,citesort}
\usepackage[psamsfonts]{eucal}
\usepackage{float}

\usepackage{graphicx,epstopdf}
\usepackage{multirow}
\usepackage{amsfonts,amsmath,amssymb,latexsym}

\usepackage{slashbox,pict2e}
\usepackage{makecell}
\usepackage{booktabs}

\usepackage{BernsteinStyle2}

\usepackage{xparse}

\newtheorem{lma}{Lemma}
\newtheorem{prp}{Proposition}
\newtheorem{theom}{Theorem}

\newcounter{Example}
\newcommand{\Example}[1]{\refstepcounter{Example}}

\DeclareMathOperator*\limif{\underline{lim}}
\DeclareMathOperator*\limsp{\overline{lim}}
\DeclareMathSizes{10}{8.95}{5.5}{5}

\begin{document}

\begin{frontmatter}

 \title{Convergence and Consistency  of\\  Recursive Least Squares with  Variable-Rate Forgetting} 

\author[Michigan]{Adam L. Bruce}\ead{admbruce@umich.edu}, 
\author[Michigan]{Ankit Goel}\ead{ankgoel@umich.edu},
and \author[Michigan]{Dennis S. Bernstein}\ead{dsbaero@umich.edu}

\address[Michigan]{Department of Aerospace Engineering, The University of Michigan, Ann Arbor, Michigan 48109, United States}



\begin{abstract}
A recursive least squares algorithm with variable rate forgetting (VRF) is derived by minimizing a quadratic cost function.
Under persistent excitation and boundedness of the forgetting factor, the minimizer given by VRF is shown to converge to the true parameters.
In addition, under persistent excitation and with noisy measurements, where the noise is uncorrelated with the regressor, conditions are given under which the minimizer given by VRF is a consistent estimator of the true parameters.
The results are illustrated by a numerical example involving abruptly changing parameters.
\end{abstract}

\end{frontmatter}

\section{Introduction}

Recursive least squares (RLS) is one of the foundational algorithms of systems and control theory, especially for signal processing, identification, and adaptive control \cite{ljung:83,astrom}.
An early exposition of RLS is given in \cite{albert1965}.

Standard RLS employs a constant forgetting factor $\lambda$, which enhances the importance of recent data over older data.
Although $\lambda$ can be set by the user, the performance of RLS is often extremely sensitive to the chosen value.
Consequently, choosing a suitable value of $\lambda$ is typically a trial and error process.

To remedy this problem, various techniques have been proposed to automatically vary the forgetting factor in response to the fit error.  
In particular, \cite{FortescueKershenbaumYdstie} reports a method for sequentially updating the forgetting factor to conserve the amount information used in the estimate, and \cite{PaleologuBenestySilviu} reports an update-based algorithm that uses noise statistics to control the forgetting factor.
\cite{LeungAndSo} gives a gradient-based algorithm for computing a forgetting factor that locally minimizes the mean-square error of the estimate, and \cite{SongLimBaekSung} derives a Newton-type gradient-descent algorithm that combines sequential estimation with minimization of the mean-squared error.
Finally, \cite{ParkJunKim} gives a formula based on exponentiation of the squared residual.

The present paper approaches the problem of varying the forgetting factor by deriving a generalization of RLS that includes time-dependent cost scaling and regularization.
This formulation involves a growing-window cost function, and thus is distinct from the formulation of \cite{AliHoaggMossbergBernstein}, which uses a sliding-window cost function.
The growing-window cost function is advantageous since it directly generalizes traditional RLS and has the ability to weigh recent data more heavily than older data.

The first contribution of the paper is given by Theorem 1, which introduces RLS with variable-rate forgetting (VRF), a novel extension of RLS in which the role of the constant forgetting factor $\lambda$ in RLS is replaced by a variable forgetting factor $\beta_k.$
By setting $\beta_k=\frac{1}{\lambda}$ for all $k,$ VRF specializes to RLS with constant-rate-forgetting (CRF).
The variable-rate-forgetting extensions of RLS given in \cite{FortescueKershenbaumYdstie,PaleologuBenestySilviu,LeungAndSo,SongLimBaekSung,ParkJunKim} are special cases of Theorem \ref{thm:genmin} with specific choices of $\beta_k$. 
In addition, Theorem 1 refines the variable-rate weighting used  in \cite[pp. 17, 18]{ljung:83}. 
In particular, we factor $\alpha_k$ in \cite[Eq. (2.12)]{ljung:83} as $\beta_k\cdots\beta_0,$ where $1/\beta_k$ serves as the instantaneous forgetting factor at step $k.$
This formulation allows the user to specify $\beta_k$ at each step based on the current residual or knowledge of system changes.
The second and third contributions of this paper are given by Theorems \ref{thm:convergence}, \ref{thm:consistency}, and Corollary \ref{cor:consistiflimsupzero}, which prove conditions on $\beta_k$ ensuring convergence under the assumption of persistency (Theorem \ref{thm:convergence}) and consistency under the assumption of persistency and that the regressor and sensor noise are uncorrelated (Theorem \ref{thm:consistency}, Corollary \ref{cor:consistiflimsupzero}). Specific examples of $\beta_k$ for consistent and non-consistent algorithms are given in Corollary \ref{cor::conditionsforconsistency}. %
The fourth contribution is two choices of $\beta_k$ that may be useful in practice. In Section \ref{sec:example}, we demonstrate these choices on an abruptly changing system with and without measurement noise and compare the performance of VRF and CRF for the given example. 
%

%
The notation used throughout this paper is as follows. The symbols $\BBS^n,$ $\BBN^n,$ and $\BBP^n$ denote the sets of real $n\times n$ symmetric, positive-semidefinite, and positive-definite matrices, respectively.  
For all $A\in \BBS^n$, $\lambda_i(A)$ denotes the $i$th largest eigenvalue of $A$, $\lambda_{\max}(A) \isdef \lambda_1(A)$, and $\lambda_{\min}(A) \isdef \lambda_n(A)$.
$\lfloor x \rfloor$ denotes the greatest integer less than or equal to $x\in \BBR$. Finally, for all $k \ge 0$ and $N \ge 0$, we define $\xi(k,N)\isdef \lfloor \frac{k}{N+1}\rfloor$. 

\section{Problem Formulation}
Let $\lambda\in(0,1]$, $\theta_0 \in \BBR^{n}$, and
$P_0 \in \BBP^n$.
Furthermore, for all $k \geq 0$, let
$\phi_k\in \BBR^{p\times n}$, $y_k \in\BBR^p$, $e_k \isdef y_k - \phi_k\theta$  
and define $J_k\mspace{-1mu}\colon\BBR^n\to[0,\infty)$ by
\begin{align}
    J_k(\theta) &\isdef   \sum_{i=0}^{k}\lambda^{k-i}\|e_k\|^2
    %
    + \lambda^{k+1}(\theta- \theta_0)^\rmT P_0^{-1} (\theta- \theta_0 ).\label{eq:rlscost}
\end{align}
Equation \eqref{eq:rlscost} is the cost function for CRF, the minimization of which produces the least squares estimate of $\theta$ given $y_0,\dots,y_k$.
 Since $J_k$ is quadratic and strictly convex, it follows that its unique global minimizer, 
%
    $\theta_{k+1} \isdef {\rm argmin}_{\theta\in\BBR^n} J_{k}(\theta),$
%
 is the only local minimizer. The following proposition gives the traditional RLS update equations  for computing $\theta_{k+1}$ 
 \cite{aseemrls,astrom,ljung:83}.

\begin{prp}\label{tradrls}
Under the notation and assumptions of the preceding paragraph, for all $k \ge 0$, define $J_k\mspace{-1mu}\colon\BBR^n\to[0,\infty)$ by \eqref{eq:rlscost}. 
Then 
\begin{align}
    \theta_{k+1} &= \theta_k + P_{k+1}\phi_k^\rmT(y_k - \phi_k\theta_k),\label{eq:crf_theta}
\end{align}
where
\begin{align}
    P_{k+1} &= \frac{1}{\lambda}P_k - \frac{1}{\lambda}P_k\phi_k^\rmT(\lambda I_p + \phi_kP_k\phi_k^\rmT)^{-1}\phi_kP_k.\label{eq:crf_P}
\end{align}
\end{prp}
In this paper, we introduce a generalization of \eqref{eq:rlscost} in which the forgetting factor is variable, prove a result analogous to Proposition \ref{tradrls} for the generalization, and analyze convergence and consistency for the family of algorithms thus obtained. To generalize \eqref{eq:rlscost}, for all $k\ge0$, let $\beta_k > 0$, define
\begin{align}\label{eq:rhodefn}
    \rho_k \isdef \prod_{i = 0}^k\beta_i, \quad \rho_{-1}\isdef 1,
\end{align}
and define the cost function $J_k\colon\BBR^n\to[0,\infty)$ by
\begin{align}
    J_k(\theta) &\isdef   \sum_{i=0}^{k}\frac{\rho_i}{\rho_k}\|e_k\|^2
    + \frac{1}{\rho_k}(\theta- \theta_0)^\rmT P_0^{-1} (\theta- \theta_0 ).\label{eq:vrfcost}
\end{align}
Since \eqref{eq:vrfcost} is quadratic and strictly convex, like  \eqref{eq:rlscost}, its unique global minimizer is the only local minimizer. Theorem \ref{thm:genmin} provides recursive update equations for this minimizer.

\section{RLS with Variable-Rate Forgetting}

Note that \eqref{eq:vrfcost} can be written as 
\begin{align}
    J_k(\theta) = \theta^\rmT A_k \theta - 2b_k^\rmT\theta + c_k,
\end{align}
 where
 \begin{align}
     A_k &\isdef \sum_{i=0}^k\frac{\rho_i}{\rho_k}\phi_i^\rmT\phi_i + \frac{1}{\rho_k}P_0^{-1},\label{eq:akdefn}\\
     b_k &\isdef \sum_{i=0}^k \frac{\rho_i}{\rho_k}\phi_i^\rmT y_i + \frac{1}{\rho_k}P_0^{-1}\theta_0,\label{eq:bkdefn}\\
     c_k &\isdef \sum_{i=0}^k\frac{\rho_i}{\rho_k} y_i^\rmT y_i + \frac{1}{\rho_k}\theta_0^\rmT P_0^{-1}\theta_0. 
 \end{align}
Since $A_k$ is positive definite, we define the positive-definite matrix
 \begin{align}
     P_{k} \isdef A_{k-1}^{-1}, \label{eq:pkdefn}
 \end{align}
where $A_{-1} \isdef P_0^{-1}$.

The following result, {\it RLS with variable-rate forgetting (VRF)}, generalizes Proposition \ref{tradrls} to the minimizer of  \eqref{eq:vrfcost}.

\begin{theom}\label{thm:genmin}
Let $\theta_0 \in \BBR^{n}$, $P_0 \in \BBP^n$,  
and, for all $k \geq 0$, let
$\phi_k\in \BBR^{p\times n}$,
$y_k \in\BBR^p,$ 
and $\beta_k\in (0,\infty)$.
Then the minimizer $\theta_{k+1}$ of \eqref{eq:vrfcost} is given by
\begin{align}
    \theta_{k+1} &= \theta_{k} + P_{k+1} \phi_k^\rmT(y_k - \phi_k \theta_{k}),\label{eq:vrf_theta}
\end{align}
and
\begin{align}
    P_{k+1} &= L_k - L_k\phi_k^\rmT(I_p  + \phi_k L_k\phi_k^\rmT)^{-1}\phi_k L_k\label{eq:vrf_P},\\
    L_k &\isdef \beta_k P_k\label{eq:vrf_L}.
\end{align}
\end{theom}
The proof of Theorem \ref{thm:genmin} requires the following lemma.
\begin{lma}\label{lem:pinvupdatevrf}
    Let $P_0 \in \BBP^n$ and, for all $k\ge0$, let $\beta_k > 0$, define $\rho_k$ by \eqref{eq:rhodefn}, and define $P_k$ by \eqref{eq:pkdefn}.  Then, for all $k\ge0,$
    \begin{align}
         P_{k+1}^{-1} &= \frac{1}{\beta_k}P_k^{-1} + \phi_k^\rmT\phi_k\label{eq:pinvupdatevrf}\\
         &= \frac{1}{\rho_k}\left(P_0^{-1} + \sum_{i=0}^k\rho_i\phi_i^\rmT\phi_i\right).\label{eq:pinvsum}
    \end{align}
\end{lma}

\textbf{Proof.} Let $k\ge0.$  It follows from \eqref{eq:akdefn} that
%
    $A_k 
    %
    %
    =\frac{1}{\beta_k}A_{k-1} + \phi_k^\rmT\phi_k,$ 
%
which, using \eqref{eq:pkdefn}, implies \eqref{eq:pinvupdatevrf}.
Furthermore, \eqref{eq:pinvupdatevrf} implies $P_1^{-1} = \frac{1}{\rho_0}(P_0^{-1} + \rho_0\phi^\rmT_0\phi_0),$ which confirms \eqref{eq:pinvsum} for $k = 0$. Next, let $k > 0$ and suppose for induction that \eqref{eq:pinvsum} holds for $k-1$. From \eqref{eq:pinvupdatevrf} it follows that
%
    $P_{k+1}^{-1} = \frac{1}{\beta_k}P_k^{-1} + \phi_k^\rmT\phi_k
    %
    %
    =  \frac{1}{\rho_k}\left(P_0^{-1} + \sum_{i=0}^{k-1}\rho_i\phi_i^\rmT\phi_i\right) + \frac{\rho_k}{\rho_k}\phi_k^\rmT\phi_k
    %
    = \frac{1}{\rho_k}\left(P_0^{-1} + \sum_{i=0}^{k}\rho_i\phi_i^\rmT\phi_i\right).$
    \hfill\mbox{$\square$}

%
%
%
%
%
%
%
%
%
%
%
\textbf{Proof of Theorem \ref{thm:genmin}.} Let $k \ge0$. To prove \eqref{eq:vrf_P}, note that it follows from \eqref{eq:vrf_L}, \eqref{eq:pinvupdatevrf}, and the matrix inversion lemma that
%
%
 %
     $P_{k+1} 
     %
     %
     %
     =\left(\frac{1}{\beta_k}P_k^{-1}+\phi_k^\rmT\phi_k\right)^{-1}
     %
     = L_k - L_k\phi_k^\rmT\left(I_p  + \phi_k L_k\phi_k^\rmT\right)^{-1}\phi_k L_k.$ 
%
To prove \eqref{eq:vrf_theta}, note that
%
%
\eqref{eq:bkdefn}, \eqref{eq:pkdefn}, and \eqref{eq:pinvupdatevrf} imply that
\begin{align}
     &\theta_{k+1} 
     %
     %
     %
     =P_{k+1}\left(\phi_k^\rmT y_k+\frac{\rho_{k-1}}{\rho_k}\left[\sum_{i=0}^{k-1} \frac{\rho_i}{\rho_{k-1}}\phi_i^\rmT y_i + \frac{1}{\rho_{k-1}}P_0^{-1}\theta_0\right]\right)  \nn\\   
     %
     %
     %
     %
     &=P_{k+1}\left(\phi_k^\rmT\phi_k + \frac{1}{\beta_k}P_k^{-1}\right)\theta_k+P_{k+1}\phi_k^\rmT (y_k-\phi_k\theta_k)\nn\\
     &=\theta_k + P_{k+1}\phi_k^\rmT (y_k-\phi_k\theta_k),\nn \tag*{\mbox{$\square$}} %
\end{align}
%
%
For all $k\ge0,$ let $\beta_k = \frac{1}{\lambda}$. 
Then \eqref{eq:vrfcost} specializes to \eqref{eq:rlscost}, and \eqref{eq:vrf_theta}--\eqref{eq:vrf_L} specialize to \eqref{eq:crf_P} and \eqref{eq:crf_theta}.  Theorem \ref{thm:genmin} thus includes Proposition \ref{tradrls} as a special case.

\section{Convergence of VRF} \label{sec:conv}

\begin{defn}
A sequence $(S_k)_{k\ge 0} \subset \BBN^n$ is {\rm{persistent}} if there exist $N \ge 1$ and $\alpha>0$ 
such that, for all $j \ge 0$,
\begin{align}
    \alpha I_n \leq \sum_{i=0}^{N}S_{i+j}.\label{eq:persistdefn}
\end{align}
%
The numbers $\alpha$ 
and $N$ are, respectively, the {\rm lower bound} 
 and {\rm persistency window} of $(S_k)_{k\ge0}$.
The sequence $(\phi_k)_{k\ge 0} \subset  \BBR^{n\times m}$ is {\rm{persistent}} if $(\phi_k^\rmT \phi_k)_{k\ge 0}$ is persistent. 
\end{defn}
\begin{theom}\label{thm:convergence} Let $(\phi_k)_{k\geq 0}\subset  \BBR^{n\times m}$,
be persistent, 
%
let $\theta\in\BBR^n$, 
and, for all $k\geq 0$, let $y_k = \phi_k\theta$.
%
Furthermore, let $a > 1$ and, for all $k \geq 0,$ 
let $ \beta_k \ge 1$. 
Finally, let $\theta_0 \in \BBR^n$, 
let $P_0 \in \BBP^{n}$, 
and, for all $k \geq 0$, define $\theta_{k+1}$ by \eqref{eq:vrf_theta}--\eqref{eq:vrf_L}.  
Then $\lim_{k\rightarrow\infty} \theta_k = \theta.$
\end{theom}

Let $k\ge0$ 
and define $\tilde\theta_k \isdef \theta_k - \theta.$
Using \eqref{eq:vrf_theta} 
and \eqref{eq:pinvupdatevrf} it follows that
%
    $\tilde\theta_{k+1} 
    %
    %
    =  (I_n- P_{k+1}\phi_k^\rmT\phi_k)\tilde\theta_k
    %
    %
    =  \frac{1}{\beta_k}P_{k+1}P_{k}^{-1}\tilde\theta_k,$ 
%
%
thus
%
    $\tilde\theta_k = \frac{1}{\rho_{k-1}}P_k P_0^{-1}\tilde\theta_0.$ 
From \eqref{eq:pinvsum}, it follows that
\begin{align}
    &\limsp_{k\to\infty}\|\tilde\theta_k\|^2 \leq \limsp_{k\to\infty}\frac{\lambda_{\max}(P_k^2)}{\rho_{k-1}^2}\|P_0^{-1}\theta_0\|^2\nn\\
    &\leq \limsp_{k\to\infty}\frac{\|P_0^{-1}\theta_0\|^2}{\lambda^2_{\max}\left(P_0^{-1}+\sum_{i=0}^{k-1}\rho_i\phi_i^\rmT\phi_i\right)} \nn\\
    &
    \leq \limsp_{k\to\infty}\frac{\|P_0^{-1}\theta_0\|^2}{[\lambda_{\max}(P_0^{-1})+\xi(k,N+1)\alpha]^2} = 0.\nn \tag*{\mbox{$\square$}}
\end{align}
\section{Consistency of VRF} \label{sec:consist}
A sequence $(X_k)_{k\ge0}$ of vector-valued random variables on $\Omega$ is a {\it consistent estimator} of $\theta\in\BBR^n$ if, for all $\varepsilon >0$, 
\begin{align}
    \lim\limits_{k\to\infty}\BBP(\{\omega\in\Omega\mspace{-1mu}\colon  \|X_k(\omega)- \theta\| < \varepsilon\}) = 1.
\end{align}
When $\theta$ is understood, for brevity, we call such sequences \textit{consistent}.
In this section we give conditions on $\beta_k$ which are necessary and sufficient for the consistency of VRF when the measurements of $\phi_k \theta$ are corrupted by noise.  
\begin{defn}
Let $(S_i)_{i \ge 0}\subset \BBN^n$ be persistent with lower bound $\alpha$ and window $N$.  Then the {\rm upper bound} $\beta\in(0,\infty)\cup\{\infty\}$ of $(S_i)_{i \ge 0}$ is
\begin{equation}
    \beta \isdef\sup_{j\ge0} \lambda_{\rm max} \left( \sum_{i=0}^N S_{i+j}\right).\label{eq:betadefn}
\end{equation}
\end{defn}
%

\begin{lma}\label{lem:sumbounds}
Let $(S_i)_{i \ge 0}\subset \BBN^n$ be persistent with  window $N$, lower bound $\alpha$, and upper bound $\beta,$ and let $(a_i)_{i\ge 0}$ be a nondecreasing sequence of nonnegative numbers. Then, for all $k \ge 0$,
\begin{align}
    \alpha \ell_{\xi(k,N)-1} I_n \leq \sum_{i=0}^{k}a_iS_i \leq \beta r_{\xi(k,N)}I_n,\label{eq:boundLemma}
\end{align}
where
%
%
%
%
    %
    $\ell_j \isdef\sum_{i=0}^{j}a_{i(N+1)}$ and $r_j \isdef \sum_{i=0}^{j}a_{i(N+1)+N}.$
    %
%
%
\end{lma}

 \textbf{Proof.} In the case where $\beta = \infty$, the upper bound of \eqref{eq:boundLemma} is immediate. Hence, assume $\beta <\infty$. Let $k \geq 0$. Since $(a_i)_{i \ge 0}$ is nondecreasing, for all $j \ge 0$ and $i \in\{ 0,\dots,N\}$, $a_{i+j} \leq a_{N+j}$ and $a_{j} \leq a_{i+j}$. From \eqref{eq:persistdefn} and \eqref{eq:betadefn} it follows that
 %
      $\alpha a_j I_n \leq a_{j}\sum_{i=0}^{N}S_{i+j} \leq \sum_{i=0}^{N}a_{i+j}S_{i+j},$  
%
%
thus
%
    %
    $\alpha \ell_{\xi(k,N)-1} I_n
    %
    %
    \leq \sum_{q=0}^{\xi(k,N)-1}\sum_{i=0}^Na_{i+q(N+1)}S_{i+q(N+1)} 
    %
    %
    \leq \sum_{i=0}^{k}a_iS_i.$ 
 %
 %
%
Similarly,      $\sum_{i=0}^{N}a_{i+j}S_{i+j} \leq a_{N+j}\sum_{i=0}^{N}S_{i+j} \leq a_{N+j}\beta I_n,$ thus
%
    $\sum_{i=0}^{k}a_iS_i 
    %
    \leq \sum_{q=0}^{\xi(k,N)-1}\mspace{-8mu}a_{q(N+1)+N}\beta I_n + a_k\beta I_n\leq r_{\xi(k,N)} \beta I_n.$ \hfill\mbox{$\square$}
%
%
%
\begin{theom}\label{thm:consistency}
Let $(\phi_k)_{k\geq 0}$ be a persistently exciting sequence with window $N$,
lower bound $\alpha$, and upper bound $\beta<\infty$. 
Let $\theta\in\BBR^n$, 
$P_0\in \BBP^{n}$, 
and $\theta_0 \sim \mathcal{N}(\theta,P_0)$. 
Let $(\nu_k)_{k\ge0}$ be an $\BBR^p$-valued stationary Gaussian white-noise process with variance $V$ and uncorrelated with $\theta_0,$ and
define $y_k = \phi_k\theta + \nu_k. $
%
%
Furthermore, for all $k \geq 0$, 
let $\beta_k \geq 1$,
and define $\theta_{k+1}$ by \eqref{eq:vrf_theta}--\eqref{eq:vrf_L}. 
%
%
Then, for all $k\ge0,$ $\theta_k$ is a Gaussian random variable with mean $\bar{\theta}$.
%
%
%
Then 
\begin{align}
    &\frac{\alpha\lambda_{\min}(V )}{\beta^2}\limif\limits_{k\to\infty}\frac{q_{\rml,\xi(k,N)}}{s_{\rmu,\xi(k,N)}^2}\label{eq:covboundlow} 
    \leq \limif\limits_{k\to\infty}\lambda_{\min}({\rm var}(\theta_k))\\
    &\leq \limsp\limits_{k\to\infty}\lambda_{\max}({\rm var}(\theta_k))
    %
    \leq \frac{\beta\lambda_{\max}(V )}{\alpha^2}\limsp\limits_{k\to\infty}\frac{q_{\rmu,\xi(k,N)}}{s_{\rml,\xi(k,N)}^2},\label{eq:covboundup}
\end{align}
where, for all $j\ge0,$ $s_{\rml,j} \isdef \sum_{i=0}^{j-1}\rho_{i(N+1)}, %
        \quad s_{\rmu,j}\isdef \sum_{i=0}^j\rho_{i(N+1)+N}$, $s_{\rmu,j}\isdef \sum_{i=0}^j\rho_{i(N+1)+N}$, $q_{\rml,j} \isdef \sum_{i=0}^{j-1}\rho^2_{i(N+1)}$, and $ q_{\rmu,j}\isdef \sum_{i=0}^j\rho^2_{i(N+1)+N}$.
\end{theom}
\textbf{Proof.}
With base case $\theta_0\sim \mathcal{N}(\theta,P_0)$, suppose for induction that $\theta_k \sim \mathcal{N}(\theta,{\rm var}(\theta_k)).$ Define $\tilde \theta_k \isdef \theta_k - \theta.$ From \eqref{eq:vrf_theta}, it follows that
%
    $\tilde\theta_{k+1} 
    %
    = \beta_k^{-1}P_{k+1}P_k^{-1}\tilde\theta_k + P_{k+1}\phi_k^\rmT \nu_k.$ 
%
Since $\theta_k \sim \mathcal{N}(\theta,{\rm var}(\theta_k)),$ it follows from Lemma \ref{lem:GCV} 
that $\tilde\theta_k \sim \mathcal{N}(0,{\rm var}(\theta_k))$. Next, define $z_k \isdef P_k^{-1}\tilde\theta_k$.
%
%
Since $\tilde\theta_k \sim \mathcal{N}(0,{\rm var}(\theta_k)),$ it follows from Lemma \ref{lem:GCV} 
that $z_k \sim \mathcal{N}(0,P_k^{-1}{\rm var}(\theta_k)P_k^{-1})$. Since $\nu_k$ is uncorrelated with $\nu_0,\ldots,\nu_{k-1},\theta_0$, it follows that $\nu_k$ and $z_k$ are also uncorrelated. Furthermore,
%
    $z_{k+1} = \beta_k^{-1}z_k + \phi_k^\rmT \nu_k,$
    %
%
and thus $[z_k  \  \nu_k]^\rmT\sim \mathcal{N}(0_{2\times 1},{\rm diag}({\rm var}(z_k),V))$.
%
%
Therefore, Lemma \ref{lem:GCV} 
implies that $z_{k+1} \sim \mathcal{N}(0,{\rm var}(z_{k+1}))$ and  
%
  ${\rm var}(z_{k+1}) = \beta_k^{-2}{\rm var}(z_k) + \phi_k^\rmT V  \phi_k.$ 
%
Since $\theta_{k+1} = P_{k+1}z_{k+1} + \theta$, it follows from Lemma \ref{lem:GCV} 
that $\theta_{k+1}\sim \mathcal{N}(\theta,P_{k+1}{\rm var}(z_{k+1})P_{k+1})$.
Thus, for all $k\ge0$, $\theta_k$ is a Gaussian random variable with mean $\theta$. 
Since ${\rm var}(z_{0}) = P_0^{-1}P_0P_0^{-1} = P_0^{-1}$, it follows that
%
    ${\rm var}(z_{k+1}) = \rho_k^{-2}\left(P_0^{-1} + \sum_{i=0}^k \rho_i^2 \phi_i^\rmT V  \phi_i\right).$ 
%
For convenience, define
%
    $M_k \isdef \sum_{i=0}^k \rho_i\phi_i^\rmT\phi_i,$ 
    $M_{\nu,k} \isdef \sum_{i=0}^k \rho_i^2 \phi_i^\rmT V  \phi_i,$ 
    %
    %
   $H_{0,k} \isdef (P_0^{-1}+M_k)^{-1}P_0^{-1}(P_0^{-1}+M_k)^{-1},$
   %
   $H_{\nu,k} \isdef (P_0^{-1}+M_k)^{-1}M_{\nu,k}(P_0^{-1}+M_k)^{-1}.$
   %
   %
%
For all $k \geq 0$, it follows from Lemma \ref{lem:sumbounds} that 
\begin{align}
    \alpha s_{\rml,\xi(k,N)}I_n &\leq M_{k} \leq \beta s_{\rmu,\xi(k,N)} I_n,\label{eq:Mkineq}\\
    \alpha\lambda_{\min}(V ) q_{\rml,\xi(k,N)}I_n &\leq M_{\nu,k} \leq \beta\lambda_{\max}(V ) q_{\rmu,\xi(k,N)} I_n.\label{eq:Mnukineq}
\end{align}
Since $\beta_k \geq 1$, it follows that $q_{\ell,\xi(k,N)}\to\infty$ as $k \to \infty$, and thus $\lambda_{\max}(M_k)\to\infty$ as $k \to \infty$. From this result and Lemma \ref{lem:eigenineqlemma} 
it follows that
%
    $\limsp_{k\to\infty}\lambda_{\max}(H_{0,k})
    %
    \leq \limsp_{k\to\infty}\lambda_{\max}(P_0^{-1})/\lambda_{\max}(M_k)^2 = 0.$ 
%
Hence, $\limsp_{k\to\infty}\lambda_{\max}(H_{0,k})=0.$
Noting that ${\rm var}(\theta_k) = H_{0,k}+H_{\nu,k}$,
it follows from Lemmas \ref{lma:limsupsubadd} and \ref{lem:eigenineqlemma}, 
\eqref{eq:Mkineq}, and 
\eqref{eq:Mnukineq} that
\begin{align}
    &\limsp\limits_{k\to\infty} \ \lambda_{\max}({\rm var}(\theta_k))
    %
    \leq \limsp\limits_{k\to\infty}\lambda_{\max}(H_{0,k}) + \limsp\limits_{k\to\infty}\lambda_{\max}(H_{\nu,k})
    \nn\\
    &=\limsp\limits_{k\to\infty}\lambda_{\max}(H_{\nu,k})
    %
    \leq \limsp\limits_{k\to\infty}\frac{\lambda_{\max}(M_{\nu,k})}{\lambda_{\max}(P_0^{-1}+M_k)^2}\nn\\
    &\leq\limsp\limits_{k\to\infty}\frac{\lambda_{\max}(M_{\nu,k})}{\lambda_{\max}(M_k)^2}
    %
    \leq \frac{\beta\lambda_{\max}(V )}{\alpha^2}\limsp\limits_{k\to\infty}\frac{\ q_{\rmu,\xi(k,N)}}{s_{\rml,\xi(k,N)}^2},\nn
    %
\end{align}
%
%
Since $\limif_{k\to\infty}\lambda_{\min}(H_{0,k}) \leq \limsp_{k\to\infty}\lambda_{\max}(H_{0,k}) = 0$, it follows that $\limif_{k\to\infty}\lambda_{\min}(H_{0,k}) = 0.$ 
Thus, from Lemmas \ref{lma:limsupsubadd}, \ref{lem:eigenineqlemma}, 
and \ref{lma:liminfsuplma}, \cite[Fact 10.4.13]{matmath}, 
\eqref{eq:Mkineq}, and 
\eqref{eq:Mnukineq}, it follows that
\begin{align}
    %
    &\frac{\alpha\lambda_{\min}(V )}{\beta^2}\limif\limits_{k\to\infty}\frac{q_{\rml,\xi(k,N)}}{s_{\rmu,\xi(k,N)}^2}
    %
    \leq\limif\limits_{k\to\infty}\frac{\lambda_{\min}(M_{\nu,k})}{\lambda_{\min}(M_k)^2}\nn\\
    &=\limif\limits_{k\to\infty}\frac{\lambda_{\min}(M_{\nu,k})}{[\lambda_{\max}(P_0^{-1})+\lambda_{\min}(M_k)]^2}
    %
    \leq\limif\limits_{k\to\infty}\frac{\lambda_{\max}(M_{\nu,k})}{\lambda_{\min}(P_0^{-1}+M_k)^2}\nn\\ 
    %
    %
    &\leq\limif\limits_{k\to\infty}[\lambda_{\min}(H_{0,k}) + \lambda_{\min}(H_{\nu,k})]
    %
    = \limif\limits_{k\to\infty}\lambda_{\min}({\rm var}(\theta_k)).
    \tag*{\hfill $\square$}
\end{align}
%
 %


%
%


\begin{cor}\label{cor:consistiflimsupzero}
    Under the notation and assumptions of Theorem \ref{thm:consistency}, consider the following statements: $i$)         
    $\limsp\limits_{k\to\infty}q_{\rmu,\xi(k,N)}/s_{\rml,\xi(k,N)}^2 = 0$, $ii$) 
    $(\theta_k)_{k\ge 0}$ is consistent, $iii$), 
    $\limif\limits_{k\to\infty}q_{\rml,\xi(k,N)}/s_{\rmu,\xi(k,N)}^2 = 0$
    %
    %
    Then ${\it i})\Longrightarrow{\it ii})\Longrightarrow{\it iii}).$
\end{cor}
\textbf{Proof.} To prove ${\it i})\Longrightarrow{\it ii})$, let
$\limsp\limits_{k\to\infty}q_{\rmu,\xi(k,N)}/s_{\rml,\xi(k,N)}^2 = 0.$
Then $\lim_{k\to\infty}\lambda_{\max}({\rm var}(\theta_k)) = 0$. Thus, from Lemma \ref{lem:GCC}, it follows that $(\theta_k)_{k\ge 0}$ is consistent. To prove ${\it ii})\Longrightarrow{\it iii})$, suppose that $(\theta_k)_{k\ge 0}$ is consistent. Then, from Lemma \ref{lem:GCC}, it follows that 
$\limif_{k\to\infty}\lambda_{\min}({\rm var}(\theta_k)) = 0$, and therefore 
%
$\limif_{k\to\infty}q_{\rml,\xi(k,N)}/s_{\rmu,\xi(k,N)}^2 = 0.$ 
%
\hfill\mbox{$\square$}
\begin{cor}\label{cor::conditionsforconsistency}
    Under the notation and assumptions of Theorem \ref{thm:consistency}, the following statements hold: $i$) assume that $\prod_{k\ge 0}\beta_k$ is finite.  Then $(\theta_k)_{k\ge 0}$ is consistent; $ii$) let $\beta_0 =1$ and for all $k > 0$, let $\beta_k = 1 + \frac{1}{k}$. Then $(\theta_k)_{k\ge 0}$ is consistent; $iii$) let $\gamma \in [1,\infty)$, and, for all $k \ge 0$, let $\beta_k = \gamma.$ Then $(\theta_k)_{k\ge 0}$ is consistent if and only if $\gamma = 1.$
    %
        
\end{cor}
\textbf{Proof.} To prove $i$), suppose that $\prod_{k\ge 0}\beta_k = \rho$ and let $\varepsilon >0$. Thus there exists         $K > 0$ such that, for all $i \geq K$, $\rho - \varepsilon < \rho_{i} < \rho + \varepsilon.$ Let $k_\varepsilon > 0$ be the smallest integer such that $\xi(k_\varepsilon,N)(N+1) \geq K$, and define $B_\varepsilon \isdef \sum_{i=0}^{\xi(k_\varepsilon,N)}\rho^2_{i(N+1)+N}$ and $C_\varepsilon \isdef \sum_{i=0}^{\xi(k_\varepsilon,N)}\rho_{i(N+1)}$. Then, for all $k > k_\varepsilon$, 
\begin{align}
    %
    \dfrac{q_{\rmu,\xi(k,N)}}{s_{\rml,\xi(k,N)}^2}
    %
    %
    \leq\frac{B_\varepsilon+(\rho + \varepsilon)^2(\xi(k,N)-\xi(k_\varepsilon,N)-1)}{\left(C_\varepsilon+(\rho - \varepsilon)(\xi(k,N)-\xi(k_\varepsilon,N)-1)\right)^2}.\label{eq:limsupcor4i} 
\end{align}
Since the limit superior of the left-hand side of \eqref{eq:limsupcor4i} is zero,  
%
it follows that $(\theta_k)_{k\ge 0}$ is consistent. 
To prove $ii$), for all $k \ge 0$, let $\beta_k = 1+\frac{1}{k}$. Then, for all $i \ge 0$,
$
    \rho_{i} = i + 1,
$
%
and thus $q_{\rmu,\xi(k,N)}$ and $s_{\rml,\xi(k,N)}^2$ are polynomials of degree three and four, respectively. Hence the limit superior is zero, and therefore $(\theta_k)_{k\ge 0}$ is consistent.
To prove $iii$), suppose that $\gamma = 1$. Then 
%
    $\limsp_{k\to\infty}q_{\rmu,\xi(k,N)}/s_{\rml,\xi(k,N)}^2 = \limsp_{k\to\infty}\xi(k,N)^{-1} = 0.$
%
Hence $(\theta_k)_{k\ge 0}$ is consistent. Conversely, suppose  $\gamma > 1$. Then, for all $i \ge 0$, $\rho_i = \gamma^{i+1}$, and thus
\begin{align}
    %
    \limif\limits_{k\to\infty}&\dfrac{q_{\rml,\xi(k,N)}}{s_{\rmu,\xi(k,N)}^2}
    %
    %
    %
    =\frac{1}{\gamma^{2N}}\frac{(1-\gamma^{(N+1)})^2}{1-\gamma^{2(N+1)}}\limif\limits_{k\to\infty}\frac{1-\gamma^{2(N+1)(\xi(k,N)+1)}}{(1-\gamma^{(N+1)(\xi(k,N)+1)})^2}
    \nn\\
    &=\frac{1}{\gamma^{2N}}\frac{\gamma^{(N+1)}-1}{\gamma^{(N+1)}+1},\nn
\end{align}
which is positive because $\gamma > 1$.  Therefore, $(\theta_k)_{k\ge 0}$ is not consistent.
\hfill\mbox{$\square$}

Corollary \ref{cor::conditionsforconsistency} shows that if $\prod_{k\ge 0}\beta_k$ converges, then VRF is consistent, but also that the converse is false, since $\prod_{k>0}1+\frac{1}{k} = \infty.$
Furthermore, CRF is consistent if and only if $\lambda=1.$
\section{Example:  Abruptly Changing Parameters}\label{sec:example}
\begin{figure}[t]
    \centering
    \includegraphics[scale = 0.42]{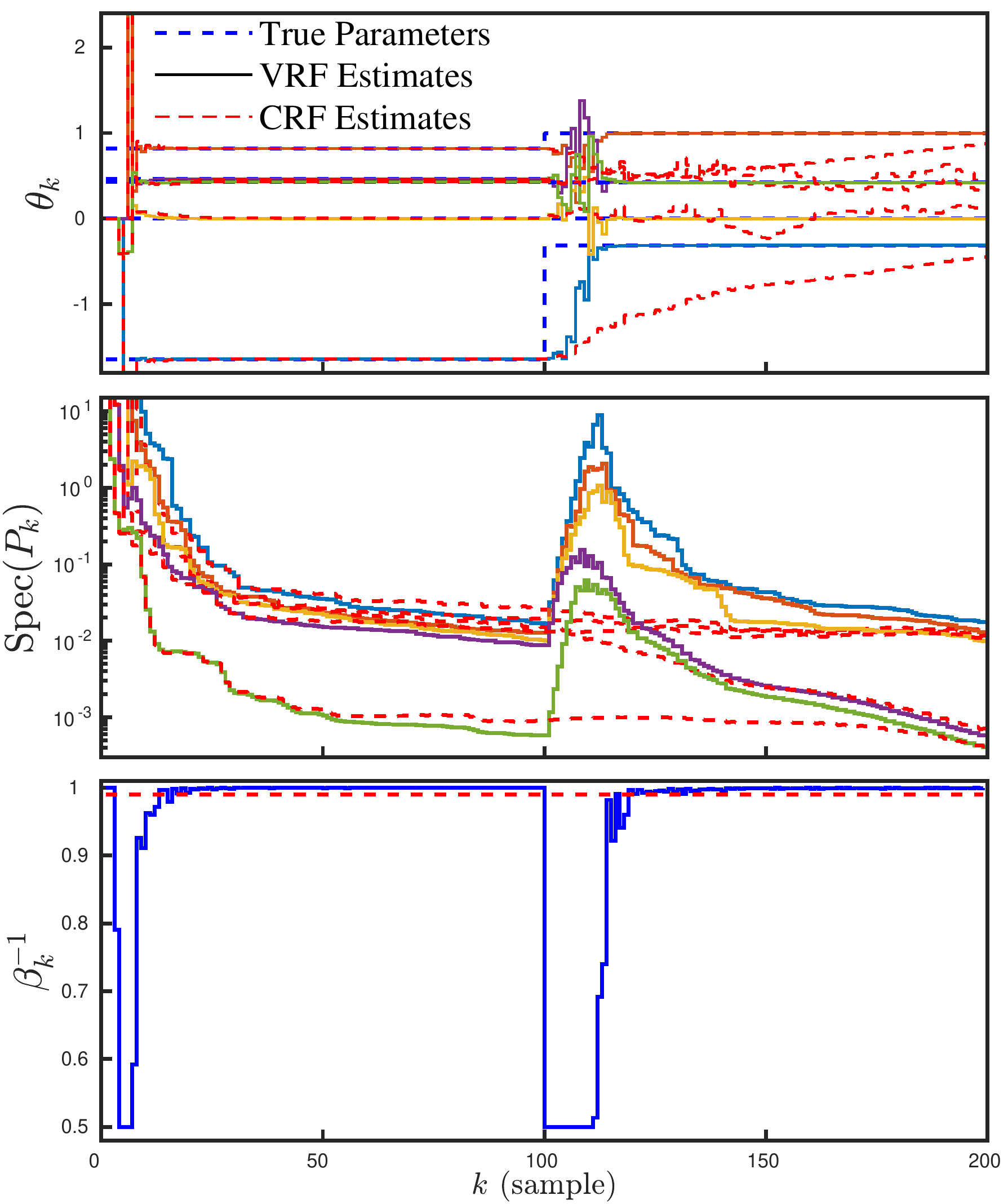}
    \vspace*{-3.5mm}
    \caption{The parameter estimate $\theta_k$ given by VRF with $\beta_k$ defined by \eqref{eq:beta_ex2} reconverges after an abrupt change in the system as guaranteed by Theorem \ref{thm:convergence}. In contrast, The parameter estimate given by CRF with $\lambda = 0.99$ requires many samples to reconverge.}
    \label{fig:example_1}
\end{figure}

\begin{figure}[t]
    \centering
    \includegraphics[scale = 0.42]{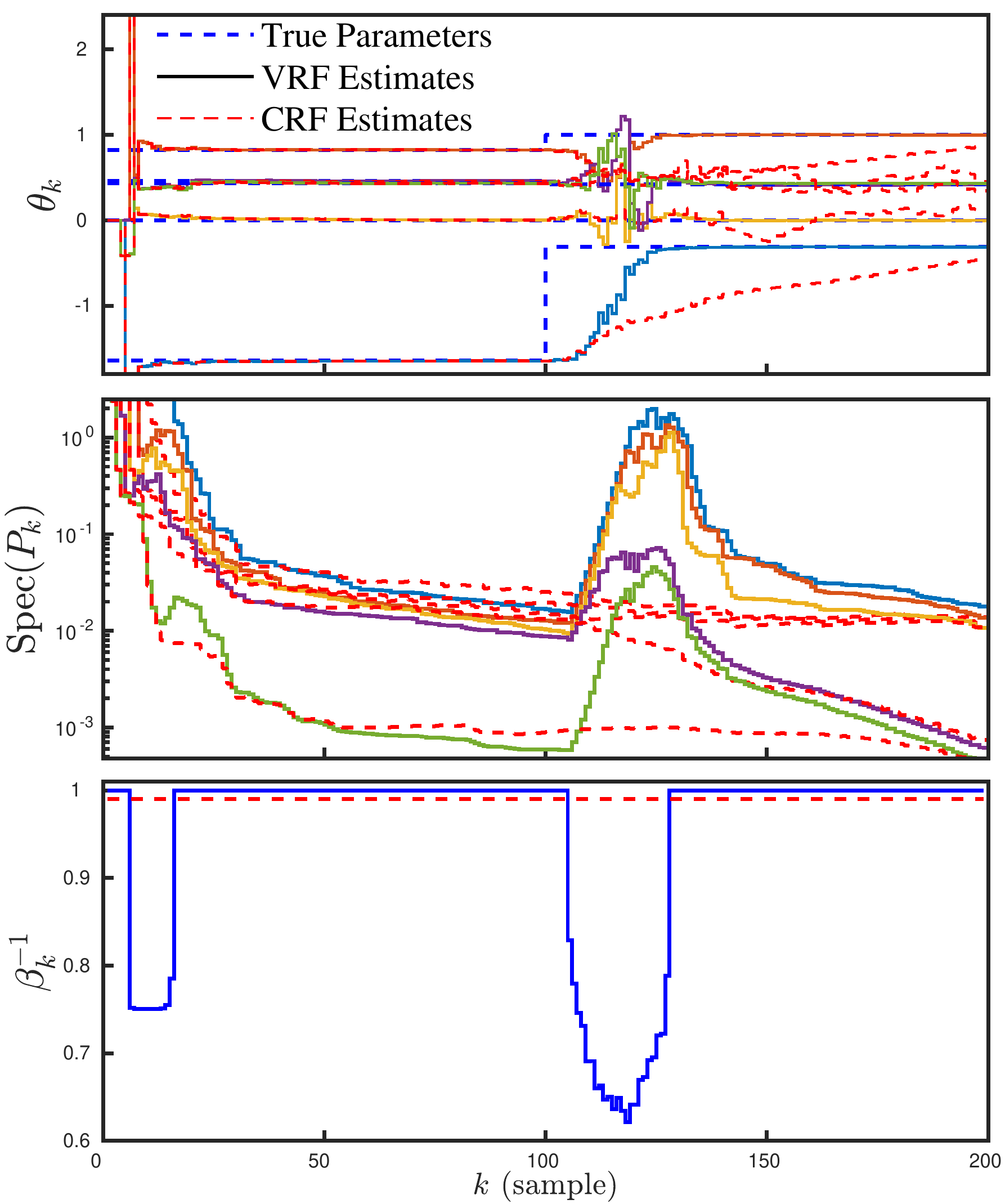}
    \vspace*{-3.5mm}
    \caption{The parameter estimate $\theta_k$ given by VRF with $\beta_k$ defined by \eqref{eq:beta_example4} reconverges after an abrupt change in the system with noisy measurements. In contrast, The parameter estimate given by CRF with $\lambda = 0.99$ requires many samples to reconverge.} 
    \label{fig:example_2}
\end{figure}

Consider a mass-spring-damper system with $m = 5$ kg, $k = 1$ N/m, and $b = 1$ N$\cdot$sec/m sampled at 1 sample/sec, and suppose that at $100$ samples the parameters of the system abruptly change to $k = 10$ N/m and $b = 0.01$ N$\cdot$sec/m. This process is modeled by the time-varying discrete-time transfer function
\begin{align}
    G_k(\mathbf{q}) = \begin{dcases}
                         \frac{0.4606 \mathbf{q} + 0.4307}{\mathbf{q}^2 - 1.64 \mathbf{q} + 0.8187}, & k < 100,\\
                         \frac{0.4218 \mathbf{q} + 0.4215}{\mathbf{q}^2 - 0.3116 \mathbf{q} + 0.998}, & k \geq 100, 
                      \end{dcases}
                       \label{eq:ex2_tf}
\end{align}
where $\mathbf{q}$ is the forward shift operator.
For all $k \ge 0$, let $u_k \sim\mathcal{N}(0,1),$
%
and define 
\begin{align}
    \beta_k \isdef 1 + \eta \ \text{sat}_\gamma(\|y_k - \phi_k\theta_k\|), \label{eq:beta_ex2}
\end{align}
where $\eta, \gamma > 0$, and $\text{sat}_\gamma$ is the unit-slope saturation function with saturation level $\gamma.$
Figure \ref{fig:example_1} shows the performance of VRF with $\gamma = \eta = 1$ and CRF with $\lambda = 0.99$. VRF converges to the initial system parameters and reconverges to the modified parameters in about 10 samples, illustrating Theorem \ref{thm:convergence}. In contrast, while CRF converges to the initial parameters, reconvergence to the modified parameters is still not achieved at $200$ samples. 
Next, consider the same system with the output corrupted by additive noise $\nu_k \sim \mathcal{N}(0,0.05)$, 
and define
\begin{align}
    \beta_k\isdef \begin{dcases}
                        1 + \eta \  \text{sat}_\gamma(E_\tau), & E_\tau > 1,\\
                        1, & E_\tau \leq 1,\\
                    \end{dcases}\label{eq:beta_example4}
\end{align}  
%
where $\tau \in \BBN$ and 
%
    $E_\tau \isdef \left(\frac{1}{\tau}\sum_{i=k-\tau}^k\|y_i-\phi_i\theta_i\|^2\right)^{1/2}.$
%
Figure \ref{fig:example_2} shows the performance of VRF with $\eta = 1$, $\gamma = 5$, and $\tau = 10$, and CRF with $\lambda = 0.99$. VRF converges to the initial parameters and then reconverges to the new parameters in roughly 30 samples. As in the previous case, CRF converges to the initial parameters, but at 200 samples has still not reconverged to the modified parameters.
\hfill \mbox{\large$\diamond$}

\section*{Acknowledgments}  This research was partially supported by AFOSR under DDDAS grant FA9550-16-1-0071
(Dynamic Data-Driven Applications Systems).

\bibliographystyle{unsrt}
\bibliography{bib_paper}


\appendix

\section{Lemmas}
\begin{lma}
\label{lem:GCV}
Let $A\in \BBR^{n\times n}$ and $b\in \BBR^n$. Let $X\sim \mathcal{N}(\mu,P)$ and define $Y \isdef AX+b.$ Then $Y\sim\mathcal{N}(A\mu+b,APA^\rmT).$
\end{lma}
%
%
%


\begin{lma}
\label{lem:GCC}
Let $(\Omega,\Sigma,P)$ be a probability space, let $\theta \in \BBR^n$, and let $(X_k\colon\Omega\to \BBR^n)_{k\ge0}$ be a sequence of random variables such that, for all $k\ge0$, $X_k\sim\mathcal{N}(\theta,\Sigma_k)$. Then $(X_k)_{k\ge0}$ is a consistent estimator for $\theta$ if and only if $\lim_{k\to\infty}\Sigma_k = 0.$
\end{lma}




\begin{lma}\label{lma:limsupsubadd}
    Let $(A_k)_{k\ge0}, (B_k)_{k\ge0} \subset (\BBN^n)$. Then
    \begin{align}
        \limsp\limits_{k\to\infty}\lambda_{\max}(A_k +B_k) \leq \limsp\limits_{k\to\infty}\lambda_{\max}(A_k) + \limsp\limits_{k\to\infty}\lambda_{\max}(B_k),\nn
        \\
        \limif\limits_{k\to\infty}\lambda_{\min}(A_k +B_k) \geq \limif\limits_{k\to\infty}\lambda_{\min}(A_k) + \limif\limits_{k\to\infty}\lambda_{\min}(B_k).\nn 
    \end{align}
\end{lma}
%

\begin{lma}\label{lem:eigenineqlemma}
    Let $A\in \BBN^n$ and $B\in\BBP^n$. Then, for all $i = 1,\dots,n$,
    \begin{align}
        \frac{\lambda_{\min}(A)}{\lambda_i(B)^2}\leq \lambda_i(B^{-1}AB^{-1}) \leq \frac{\lambda_{\max}(A)}{\lambda_i(B)^2}.\label{eq:eigenineqLemma}
    \end{align}
    Now assume that $A \in \BBP^n$. Then there exist $0 < b_1 \leq b_2$ and $0 < a_1 \leq a_2$ such that 
    \begin{align}
        a_1I_n &\leq A \leq a_2I_n,
        \label{eq:basiceigenineq1}\\
        b_1I_n &\leq B \leq b_2I_n.
        \label{eq:basiceigenineq2}
    \end{align}
    Furthermore, for all $a_1,a_2,b_1,b_2$ satisfying \eqref{eq:basiceigenineq1}, \eqref{eq:basiceigenineq2}, 
    \begin{align}
        \frac{a_1}{b_2^2}I_n &\leq B^{-1}AB^{-1} \leq \frac{a_2}{b_1^2}I_n.\label{eq:boundeigenineq} 
    \end{align}
\end{lma}

\begin{lma}\label{lma:liminfsuplma}
    Let $a \in [0,\infty)$, let $(b_k)_{k\ge 0}, (c_k)_{k\ge 0} \subset [0,\infty)$, and assume that $\lim_{k\to\infty}b_k = \infty$.  Then, for all $p \ge 0$,
    \begin{align}
        %
        \limif\limits_{k\to\infty}\frac{c_k}{(a+b_k)^p} = \limif\limits_{k\to\infty}\frac{c_k}{b_k^p}\label{eq:liminfsupeq2}. 
    \end{align}
\end{lma}
%

\end{document}